\documentclass[12pt, draftcls, onecolumn]{IEEEtran}
\IEEEoverridecommandlockouts
\usepackage[ieee]{jphmacros2e}
\usepackage{kbordermatrix}
\usepackage{algorithmic}
\usepackage{algorithm}
\usepackage{multicol} 
\usepackage{pdfsync,srcltx}
\usepackage{amsmath}
\usepackage{amsfonts}
\usepackage{amssymb}
\usepackage{graphicx,subfigure,psfrag}
\usepackage{url,ifthen}
\usepackage{cite}
\usepackage[dvips,colorlinks,bookmarksopen,bookmarksnumbered,citecolor=red,urlcolor=red]{hyperref}
\usepackage[norefs,nocites]{refcheck}%options:norefs,nocites,msgs,chkunlbld
\makeatletter\let\NAT@parse\undefined\makeatother % enbl natbib with IEEE cls

\newcommand{\R}{\mathbb{R}}

\newcommand{\N}{\mathcal{N}}
\newcommand{\Z}{\mathbf{Z}}
\renewcommand{\L}{\mathcal{L}}
\newcommand{\G}{\mathbf{G}}
\newcommand{\V}{\mathbf{V}}
\renewcommand{\E}{\mathbf{E}}
\newcommand{\scr}{\mathcal}

\allowdisplaybreaks

\begin{document}
\title{\LARGE Improving Convergence Rate of Distributed Consensus Through
  Asymmetric Weights}

\author{He Hao,
        Prabir Barooah
	\thanks{He Hao and Prabir Barooah are with Department of Mechanical and Aerospace Engineering, University of Florida, Gainesville, FL 32611, USA, email: {\tt\small {hehao,pbarooah}@ufl.edu}. This work was supported by the National Science Foundation through Grant CNS-0931885 and ECCS-0925534.}}

% make the title area
\maketitle

\begin{abstract}
  We propose a weight design method to increase the convergence rate of
  distributed consensus. Prior work has focused on symmetric weight
  design due to computational tractability. We show that with proper
  choice of \emph{asymmetric} weights, the convergence rate can be
  improved significantly over even the symmetric optimal design. In
  particular, we prove that the convergence rate in
  a lattice graph can be made independent of the size of the graph
  with asymmetric weights. We then use a
  Sturm-Liouville operator to approximate the graph Laplacian of more
  general graphs. A general weight design method is proposed based on
  this continuum approximation. Numerical computations show that the
  resulting convergence rate with asymmetric weight design is improved
  considerably over that with symmetric optimal weights and
  Metropolis-Hastings weights.
\end{abstract}

\section{Introduction}
In distributed consensus, each agent in a network updates
its state by aggregating the information from its neighbors so that
all the agents' states reach a common value. Distributed consensus
has been widely studied in recent times due to its wide ranging
applications such as multi-vehicle rendezvous, data fusion in large
sensor network, coordinated control of multi-agent system and
formation flight of unmanned vehicles and clustered satellites,
etc. (see~\cite{tsitsiklis,jadbabaie_coordination, olfati,
  ren,martinez2007motion} and references therein).

% limitation of distributed consensus- convergence rate loss due to asymmetry
The topic of this paper is the convergence rate of distributed
consensus protocols in graphs with fixed (time invariant)
topology. The convergence rate is extremely important; it determines
practical applicability of the protocol. If the convergence rate is
small, it will take many iterations before the states of all agents
are sufficiently close. Compared to the vast literature on design of
consensus protocols, however, the literature on convergence rate
analysis is meager. Convergence rate of distributed consensus in
time-varying graphs have been studied in~\cite{Olshevsky,chan2011fast,
  cao2005lower}. The related problem of mixing time of Markov chains
is studied in~\cite{boyd2004fastest}. In~\cite{LovisariARC2011},
convergence rate for a specific class of graphs, that we call L-Z
geometric graphs, are established as a function of the number of
agents. In general, the convergence rate of consensus algorithms tend
to be slow, and decreases as the number of agents increases. It is
shown in~\cite{ultrafast} that the convergence rate can be arbitrarily
fast in small-world networks. However, networks in which communication
is only possible between agents that are close enough are not likely
to be small-world.

One of the seminal works on this subject is convex
optimization of weights on edges of the graph to maximize the
consensus convergence rate~\cite{xiao2004fast,convex}. Convex
optimization imposes the constraint that the weights of the graph must
be symmetric, which means any two neighboring agents put equal weight
on the information received from each other. The convergence rate of
consensus protocols on graphs with symmetric weights degrades
considerably as the number of agents in the network
increases. In a
D-dimensional  lattice, for instance, the convergence rate is
$O(1/N^{2/D})$ if the weights are symmetric, where $N$ is the number
of agents. This result follows as a
special case of the results in~\cite{LovisariARC2011}.  Thus, the
convergence rate becomes arbitrarily small if the size of the network
grows without bound.

% comparisons with convex optimization, finite consensus
In~\cite{ren_finite,wang2010finite,sundaram2007finite}, finite-time distributed
consensus protocols are proposed to improve the performance over
asymptotic consensus. However, in general, the finite  time needed to achieve consensus depends
the number of agents in the network. Thus, for large size of networks,
although consensus can be achieved in finite time,  the time needed
to reach consensus becomes large. 

% Contribution: spectral analysis of digraph, size-independent convergence rate, continuum approximation, distributed computed weight 
In this paper, we study the problem of how to increase the convergence
rate of consensus protocols by designing \emph{asymmetric} weights on
edges. We first consider lattice graphs and derive precise formulae
for convergence rate in these graphs. In particular, we show that in
lattice graphs, with proper choice of asymmetric weights, the
convergence rate of distributed consensus can be bounded away from
zero uniformly in $N$. Thus, the proposed asymmetric design makes
distributed consensus highly scalable. In addition, we
provide exact formulae for asymptotic steady-state consensus value. With
asymmetric weights, the consensus value in general is not the average
of the initial conditions.

We next propose a weight design scheme for arbitrary 2-dimensional
geometric graphs, i.e., graphs consisting of nodes in $\R^2$. Here we use
the idea of continuum approximation to extend the asymmetric design
from lattices to geometric graphs. We show how a Sturm-Liouville
operator can be used to approximate the graph Laplacian in the case of
lattices. The spectrum of the Laplacian and the convergence rate of
consensus protocols are intimately related. The discrete weights in
lattices can be seen as samples of a continuous weight function that
appears in the S-L operator. Based on this analogy, a weight design
algorithm is proposed in which a node $i$ chooses the weight on the
edge to a neighbor $j$ depending on the relative angle between $i$ and
$j$.  Numerical simulations show that the convergence rate with
asymmetric designed weights in large graphs is an order of magnitude
higher than that with (i) optimal symmetric weights, which
are obtained by convex optimization~\cite{xiao2004fast,convex}, and
(ii) asymmetric weights obtained by Metropolis-Hastings method, which assigns
weights uniformly to each edge connecting itself to its neighbor. The
proposed weight design method is decentralized, every node can obtain
its own weight based on the angular position measurements with its
neighbors. In addition, it is computationally much cheaper than obtaining the optimal symmetric weights using convex optimization method. The proposed weight design
method can be extended to geometric graphs in $\R^D$, but in this
paper we limit ourselves to $\R^2$.

The rest of this paper is organized as follows. Section~\ref{sec:stat} presents the problem statement. Results on size-independent convergence rate on  lattice graphs with asymmetric weight  are stated in Section~\ref{sec:lattice}. Asymmetric weight design method for more general graphs appear in Section~\ref{sec:geo}. The paper ends with conclusions and future work in Section~\ref{sec:conc}.

%%%%%%%%%%%%%%%%%%%%%%%%%%%%%%%%%%%%%%%%%%%%%%%%%%%%%%%%%%%%%%%%%%%%%%%%%
\section{Problem statement}\label{sec:stat}
%%%%%%%%%%%%%%%%%%%%%%%%%%%%%%%%%%%%%%%%%%%%%%%%%%%%%%%%%%%%%%%%%%%%%%%%%
To study the problem of distributed linear consensus in networks, we first
introduce some terminologies. The network of $N$ agents is modeled by a graph $\G
= (\V,\E)$ with vertex set $\V=\{1,\dots,N\}$ and edge set $\E \subset
\V \times \V$. We use $(i,j)$ to represent a directed edge from $i$
to $j$.  A node $i$ can receive information from $j$ if and only
if $(i,j) \in \E$.  In this paper, we assume that communication is bidirectional, i.e. $(i,j) \in \E$ if and
only if $(j,i) \in \E$. For each edge $(i,j) \in \E$ in the graph, we  associate a
weight $W_{i,j} >0$ to it. The set of neighbors of $i$ is defined as
$\scr{N}_i \eqdef \{j \in \V: (i, j)\in \E\}$. The Laplacian matrix $L$ of an arbitrary graph $\G$ with edge weights $W_{i,j}$ is defined as
\begin{align}\label{eq:laplacian}
	L_{i,j}=
	\begin{cases}  	-W_{i,j} & i \neq j, (i,j) \in \E ,\\
		\sum_{k=1}^{N} W_{i,k} & i=j, (i,k) \in \E,\\
			0&\text{ otherwise.}
			\end{cases}
\end{align}

A linear consensus protocol is an iterative update law:
\begin{align}\label{eq:dynamics_individual}
	x_i(k+1) = W_{i,i} \ x_{i}(k)+\sum_{j\in \N_i} W_{i,j} \
        x_{j}(k), \quad i \in \V,
\end{align}
with initial conditions $x_i(0) \in \R$, where $k=\{0,1,2,\cdots\}$ is
the discrete time index. Following standard practice we assume the weight matrix $W$ is a stochastic matrix, i.e. $W_{i,j}\geq 0$ and $W\mathbf{1} =\mathbf{1}$, where $\mathbf{1}$ is a vector with all entries of 1. The distributed consensus protocol~\eqref{eq:dynamics_individual} can be written in the following compact form:
\begin{align}\label{eq:dynamics}
	x(k+1)=W x(k),
\end{align}
where $x(k)=[x_1(k),x_2(k), \cdots, x_N(k)]^T$ is  the
states of the $N$ agents at time $k$. It's straightforward to obtain
the following relation $L=I-W$, where $I$ is the $N\times N$ identity
matrix and $L$ is the Laplacian matrix associated with the graph with
$W_{i,j}$ as its weights on the directed edge $(i,j)$. In addition,
their spectra are related by $\sigma(L)=1- \sigma(W)$,
i.e. $\mu_\ell(L)=1-\lambda_\ell(W)$, where $\ell \in \{1,2,\cdots,
N\}$ and $\mu_\ell, \lambda_\ell$ are the eigenvalues of $L$ and $W$
respectively. The linear distributed consensus
protocol~\eqref{eq:dynamics} implies $x(k) =W^k x(0)$.  We assume $W$
is strong connected (irreducible) and primitive. In that case the
spectral radius of $W$ is $1$ and there is exactly one eigenvalue on the unit disk. Let $\pi \in \R^{1\times N}$ be the left Perron vector of $W$ corresponding to the eigenvalue of $1$, i.e. $\pi W=\pi$, $\pi_i>0$ and $\sum_{i=1}^N \pi_i=1$, we have 
\begin{align}
	\lim_{k \to \infty} W^k=\mathbf{1} \pi,
\end{align}
Therefore,  all the states of the $N$ agents asymptotically converge to a steady state value $\bar{x}$ as $k \to \infty$,
\begin{align}\label{eq:asy_value}
	\lim_{k \to \infty} x(k)=\mathbf{1} \pi x(0)= \mathbf{1} \bar{x},
\end{align}
where $\bar{x}=\sum_{i=1}^N \pi_i x_i(0)$.

One of the most important feature of linear distributed consensus is the rate of convergence to its steady state value.  It's well known that for a primitive stochastic matrix, the rate of convergence $R$ can be measured by the spectral gap $R=1-\rho(W)$, where $\rho(W)$ is the essential spectral radius of $W$, which is defined as
\begin{align*}
	\rho(W)\eqdef \max\{|\lambda|: \lambda \in \sigma(W)\setminus \{1\}\}.	
\end{align*}
 If the eigenvalues of $W$ are real and they are ordered in a non-increasing fashion such that $1=\lambda_1 \geq \lambda_2 \geq \cdots \geq \lambda_N$, then the convergence rate of $W$ is given by
\begin{align}\label{eq:esr_bound}
	R=1-\rho(W)=\min\{1-\lambda_2, 1+\lambda_N \}.
\end{align}
In addition, from Gerschgorin circle theorem, we have that $\lambda_N \geq -1+2\max_i W_{ii}$. If $\max_i W_{ii} \neq 0$, then $1+\lambda_N$ is a constant bounded
away from $0$.  Therefore, the key to find a lower bound for the convergence rate of $W$ is to find an upper bound on the second largest eigenvalue $\lambda_2$ of $W$. Equivalently, we can
find a lower bound of the second smallest eigenvalue $\mu_2$ of the
associated Laplacian matrix $L$, since $\mu_2=1-\lambda_2$.

\begin{definition}
	We say a graph $\G$ has \emph{symmetric}  weights if
        $W_{i,j}=W_{j,i}$ for each pair of neighboring agents $(i,j)
        \in \E$. Otherwise, the weights are called
        \emph{asymmetric}. \frqed
\end{definition}
If the weights are symmetric, the matrix $W$ is doubly stochastic,
meaning that each row and column sum is $1$.

The following theorem summaries the results  in~\cite{LovisariARC2011} on
the convergence rate of consensus with symmetric weights in a broad
class of graphs that include lattices.  A $D$-dimensional lattice, specifically a $N_1 \times N_2 \times \dots \times N_D$ lattice, is a graph with $N=N_1 \times N_2 \times \dots \times N_D$ nodes, in which the nodes are placed at the integer unit coordinate points of the $D$-dimensional Euclidean space and each node connects to other nodes that are exactly one unit away from it. A $D$-dimensional lattice is drawn in $\R^D$ with a Cartesian reference frame whose axes are denoted by $x_1,x_2,\cdots,x_D$. We call a graph is a \emph{L-Z geometric graph} if it can be seen as a perturbation of regular
lattice in $D$-dimensional space; each node connects other nodes
within a certain range. The formal definition is given in~\cite{LovisariARC2011}.

\begin{theorem}\label{thm:summary}
	Let $\G$ be a $D$-dimensional connected L-Z geometric graph or
        lattice and let $W$ be any doubly stochastic matrix compatible with $\G$. Then
	\begin{align}
		\frac{c_1}{N^{2/D}} \leq R \leq \frac{c_2}{N^{2/D}},
	\end{align}
	where $N$ is the number of nodes in the graph $\G$ and $c_1, c_2$ are some constants independent of $N$.  \frqed
\end{theorem}
The above theorem states that for any connected L-Z
geometric/lattice graph  $\G$, the convergence rate of
consensus with symmetric weights cannot be bounded away from $0$ uniformly with the size $N$ of the
graph. The convergence rate of the network becomes
arbitrarily slow as $N$ increases without bound. The loss of
convergence rate with symmetric information graph has also been
observed in the vehicular formation~\cite{HH_PB_CDC:10,
  PB_PM_JH_TAC:09}. In fact, another important conclusion of the
result above is that heterogeneity in weights among nodes, as long as
$W$ is symmetric, does not change the asymptotic scaling of the
convergence rate. At best it can change the constant in front of the
scaling formula (see~\cite{boyd2004fastest} also). Therefore, even centralized weight optimization scheme
proposed in~\cite{xiao2004fast,convex} - that constrain the eights to be symmetric in
order to make the optimization problem convex - will suffer from the same issue as that of
un-optimized weights on the edges. Namely, the convergence rate will decay as
$O(1/N^{2/D})$ in a $D$-dimensional lattice/L-Z geometric graph even
with the optimized weights. In the rest of the paper, we study the problem of speeding up the convergence rate by designing \emph{asymmetric} weights.

%%%%%%%%%%%%%%%%%%%%%%%%%%%%%%%%%%%%%%%%%%%%%%%%%%%%%%%%%%%%%%%%%%%%%%%%%
\section{Fast consensus on D-Dimensional lattices}\label{sec:lattice}
%%%%%%%%%%%%%%%%%%%%%%%%%%%%%%%%%%%%%%%%%%%%%%%%%%%%%%%%%%%%%%%%%%%%%%%%%
First we establish technical results (whose proofs are provided in the appendix) on the spectrum and Perron
vectors of $D$-dimensional lattices with asymmetric weights on the edges. We then
summarize their design implications at the end of section~\ref{sec:analysis-lattice}.

%%%%%%%%%%%%%%%%%%%%%%%%%%%%%%%%%%%%%%%%%%%%%%%%%%%%%%%%%%%%%%%%%%%%%%%%%
\subsection{Asymmetric weights in lattices}\label{sec:analysis-lattice}
\begin{figure}
	  \psfrag{3}{\small$3$}
	  \psfrag{1}{\small$1$}
	  \psfrag{2}{\small$2$}
	  \psfrag{N-1}{\scriptsize $N-1$}
	  \psfrag{N}{\scriptsize $N$}
	  \psfrag{e1}{\scriptsize$W_{1,2}$}
	  \psfrag{e0}{\scriptsize$W_{2,1}$}
	  \psfrag{e3}{\scriptsize$W_{2,3}$}
          \psfrag{e2}{\scriptsize$W_{3,2}$}
	  \psfrag{e5}{\scriptsize$W_{N-1,N}$}
          \psfrag{e4}{\scriptsize$W_{N,N-1}$}
	  \psfrag{x1}{$x_1$}
	  \psfrag{o}{$o$}
\centering
\includegraphics[scale = 0.3]{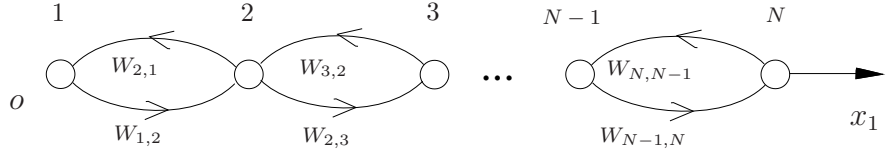}
\caption{Information graph for a 1-D lattice of $N$ agents. }
\label{fig:information_graph}

\end{figure}

We first consider distributed consensus on a 1-dimensional
lattice. This will be useful in generalizing to $D$-dimensional
lattices. Each agent interacts with its nearest neighbors in the lattice (one on each side). Its information graph is depicted in Figure~\ref{fig:information_graph}. The updating law of  agent $i$ is given by
\begin{align*}
	x_i(k+1) =W_{i,i} x_{i}(k)+ W_{i,i-1} x_{i-1}(k)+ W_{i,i+1}  x_{i+1}(k).
\end{align*}
where $i \in \{2, 3,\cdots, N-1\}$. The updating laws of the $1$-st and $N$-th agents are slightly different from the above equation, since they only have one neighbor.

The weight matrix $W^{(1)}$ for the $1$-dimensional lattice is tridiagonal:
\begin{align*}
	W^{(1)}=\begin{bmatrix}
		W_{1,1} & W_{1,2} &  \\
		W_{2,1} & W_{2,2} & W_{2,3} &\\
			&\ddots	&\ddots	&	\\
		&  W_{N-1,N-2} & W_{N-1,N-1} & W_{N-1,N} \\
		&   & W_{N,N-1} & W_{N,N} \\
	\end{bmatrix}.
\end{align*}
The following lemma gives the spectrum and the left-hand Perron vector for the weight matrix $W^{(1)}$. 

\begin{lemma}\label{lem:eig_1}
	Let $W^{(1)}$ be the weight matrix associated with the $1$-dimensional lattice with the weights given by $W_{i,i+1}=c, W_{i+1,i}=a$, where $a \neq c$ are positive constants and $a+c \leq 1$. Then its eigenvalue are 
\begin{align*}
	\lambda_1=1, \quad \lambda_{\ell}=1-a-c+2\sqrt{ac} \cos \frac{(\ell-1) \pi}{N},
\end{align*}
where $\ell \in \{2, \cdots, N\}$, and its left Perron vector is 
\begin{align*}
	\pi = 	\frac{1-c/a}{1-(c/a)^N} [1, c/a, (c/a)^2, \cdots, (c/a)^{N-1}]. \tag*{\frqed}
\end{align*}
\end{lemma}

\smallskip

\begin{figure}
\begin{center}
\psfrag{x}{\scriptsize$\ X$}
\psfrag{x1}{\scriptsize $x_1$}
\psfrag{x2}{\scriptsize $x_2$}
\psfrag{o}{\scriptsize $o$}
\psfrag{1}{\scriptsize $1$} \psfrag{2}{\scriptsize $2$} \psfrag{3}{\scriptsize $3$} 
\psfrag{4}{\scriptsize $4$} \psfrag{5}{\scriptsize $5$} \psfrag{6}{\scriptsize $6$} 
\psfrag{w21}{\scriptsize$a_1$} \psfrag{w12}{\scriptsize$c_1$} 
\psfrag{w23}{\scriptsize$c_1$} \psfrag{w32}{\scriptsize$a_1$} 
\psfrag{w41}{\scriptsize$a_2$} \psfrag{w14}{\scriptsize$c_2$} 

\includegraphics[scale = 0.5]{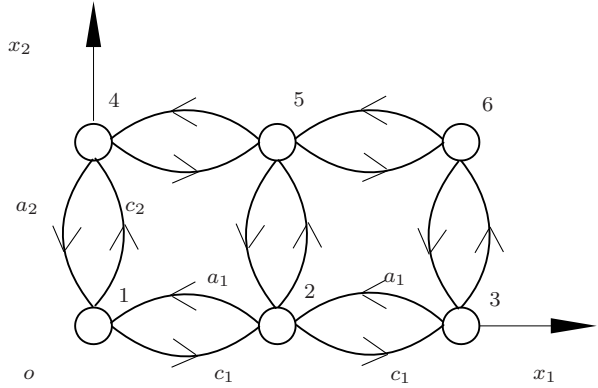}\end{center}
\caption{A pictorial representation of  a 2-dimensional lattice information graph with the weights $W^{(2)}_{i,i^{d+}}=c_d, W^{(2)}_{i,i^{d-}}=a_d$, where $d=1,2$. }\label{fig:2D}

\end{figure}

%%%%%%%%%%%%%%%%%%%%%%%%%%%%%%%%%%%%%%%%%%%%%%%%%%%%%%%%%%%%%%%%%%%%%%%%%
%\subsection{Consensus on D-dimensional lattices}\label{sec:lattice_D}
We next consider consensus on a $D$-dimensional lattice with the following weights
\begin{align}\label{eq:weights}
		W_{i,i^{d+}}=c_d,\quad W_{i,i^{d-}}=a_d, 
\end{align}
where $a_d \neq c_d$ are positive constants and $\sum_{d=1}^D a_d+c_d \leq 1$. The notation
$i^{d+}$ denotes the neighbor on the positive $x_d$ axis of node $i$
and $i^{d-}$ denotes the neighbor on the negative $x_d$ axis of node
$i$. For example, $2^{1+}$ and $2^{1-}$ in Figure~\ref{fig:2D} denote
node $3$ and node $1$, respectively, and $2^{2+}$ is node $5$.

\begin{lemma}\label{lem:eig_D}
	Let $W^{(D)}$ be the weight matrix associated with the $D$-dimensional lattice with the weights given in~\eqref{eq:weights}. Then its eigenvalues are given by
\begin{align*}
	\lambda_{\vec{\ell}} \ (W^{(D)})=1-\sum_{d=1}^{D} (1-\lambda_{\ell_d} (W^{(1)}_d)),
\end{align*}
where $\vec{\ell}=(\ell_1, \ell_2, \cdots, \ell_D)$, in which $\ell_d \in \{1,2, \cdots, N_d\}$ and $W^{(1)}_d$ is the $N_d \times N_d$ weight matrix associated with a $1$-dimensional lattice with the weights given by $W^{(1)}_d(i,i+1)=c_d,  W^{(1)}_d(i+1,i)=a_d$ and $i \in \{1,\cdots, N_d-1\}$. Its left Perron vector  is $\pi= \pi^{(1)}_D\otimes \pi^{(1)}_{D-1} \otimes \cdots \otimes \pi^{(1)}_1$, where $\pi^{(1)}_d$ is the left Perron vector of $W^{(1)}_d$.
\frqed
\end{lemma}

The next theorem shows the implications of the preceding technical results on the convergence rate in $D$-dimensional lattices. 
\begin{theorem}\label{thm:bound}
	Let $\G$ be a $D$-dimensional  lattice graph and let $W^{(D)}$ be an asymmetric stochastic matrix compatible with $\G$ with the  weights given in~\eqref{eq:weights}. Then the convergence rate satisfies	
	\begin{align}
	R \geq c_0,
	\end{align}
	where $c_0 \in (0,1)$ is a constant independent of $N$.  \frqed
\end{theorem}

\begin{remark}
  Recall from Theorem~\ref{thm:summary},  for any L-Z geometric or lattice graphs,  as long as
  the weight matrix $W$ is symmetric, no matter how do we design
  the weights $W_{i,j}$, the convergence rate
  becomes progressively smaller as the number of agents $N$ increases,
  and it cannot be uniformly bounded away from
  $0$. In contrast, Theorem~\ref{thm:bound} shows that for lattice graphs, 
  asymmetry in the weights makes the convergence rate uniformly
  bounded away from $0$. In fact, any amount of asymmetry along the
  coordinate axes of the lattice ($a_d \neq c_d$), will make this
  happen.  Asymmetric weights thus make the linear distributed consensus law highly
  scalable. It eliminates the problem of degeneration of convergence
  rate with increasing $N$. 

 The second question is where do the node states converge to with
 asymmetric weights?  Recall that the asymptotic steady
  state value of all agents is $\bar{x}=\sum_{i=1}^N \pi_i x_i(0)$. For a lattice graph, its Perron vector $\pi$ is given in Lemma~\ref{lem:eig_1} and Lemma~\ref{lem:eig_D}. Thus we can determine the  steady state value  $\bar{x}$ if the initial value $x(0)$ is given. This information is particularly useful to find the rendezvous position in multi-vehicle rendezvous problem. On the other hand, we see from Lemma~\ref{lem:eig_1} and Lemma~\ref{lem:eig_D} that if $a_d \neq c_d$, then $\pi_i \neq \frac{1}{N}$, which implies the steady-state value is not the average of the initial values. The asymmetric weight design is not applicable to distributed averaging problem.  \frqed
\end{remark}

%%%%%%%%%%%%%%%%%%%%%%%%%%%%%%%%%%%%%%%%%%%%%%%%%%%%%%%%%%%%%%%%%%%%%
\subsection{Numerical comparison}
\begin{figure}
	  \psfrag{mu}{\small $R$}
	  \psfrag{N}{\small $N$}
	  \psfrag{Convex Optimization}{\scriptsize Symmetric optimal}
	  \psfrag{Asymmetric Desgin}{\scriptsize Asymmetric design}
	  \psfrag{Prediction}{\scriptsize Lower bound~\eqref{eq:prediction}}

\centering
\includegraphics[scale = 0.35]{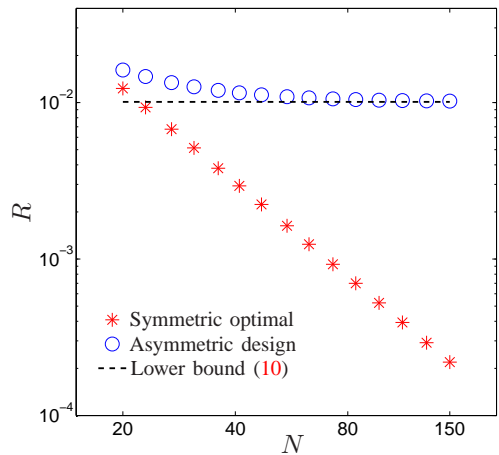}
\caption{Comparison of convergence rate of $1$-D lattice between asymmetric design and convex optimization (symmetric optimal).}
\label{fig:1Dlattice}
\vspace{-5 mm}
\end{figure}

In this section, we present the numerical comparison of the convergence rates of the distributed protocol~\eqref{eq:dynamics} between asymmetric designed weights (Theorem~\ref{thm:bound}) and symmetric optimal weights  obtained from convex optimization~\cite{xiao2004fast,convex}.   For simplicity, we take the $1$-D lattice as an example. The asymmetric weights used are  $W_{i,i+1}=c=0.3,  W_{i+1,i}=a=0.2$.  We see from Figure~\ref{fig:1Dlattice} that the convergence rate with asymmetric designed weights is much larger than that with symmetric optimal weights. In addition, given the asymmetric weight values $c=0.3, a=0.2$, we obtain  from Lemma~\ref{lem:eig_1} that $\lambda_2 \leq 0.5+2\sqrt{0.06}, \lambda_N \geq 0.5+2\sqrt{0.06}$, which implies
\begin{align} \label{eq:prediction}
R=\min\{1-\lambda_2, 1+\lambda_N \} \geq   0.5-2\sqrt{0.06}. 
\end{align}
We see from Figure~\ref{fig:1Dlattice} that the convergence rate $R$ is indeed uniformly bounded below by~\eqref{eq:prediction}. 

%%%%%%%%%%%%%%%%%%%%%%%%%%%%%%%%%%%%%%%%%%%%%%%%%%%%%%%%%%%%%%%%%%%%%
\section{Fast consensus in more general graphs}\label{sec:geo}
%%%%%%%%%%%%%%%%%%%%%%%%%%%%%%%%%%%%%%%%%%%%%%%%%%%%%%%%%%%%%%%%%%%%%
In this section, we study how to design the weight matrix $W$ to
increase the convergence rate of consensus in graphs that are more
general than lattices. We use the idea of continuum
approximation. Under some ``niceness'' properties, a graph can be
thought of as approximation of a $D$-dimensional lattice, and by
extension, of the Euclidean space corresponding to
$\R^D$~\cite{DoyleSnell}. These properties have to do with the graph
not having arbitrarily large holes etc. Precise conditions under which
a graph can be approximated by the $D$-dimensional lattice are
explored in~\cite{PB_JH_TIT:09} (for infinite graphs) and
in~\cite{LovisariARC2011} (for finite graphs). The dimension $D$ of the
corresponding lattice/Euclidean space is also determined by these
properties.

The key is to embed the discrete graph problem into a continuum-domain problem. 
We use a Sturm-Liouville operator to approximate the Laplacian matrix
of a $D$-dimensional geometric graph. A D-dimensional geometric graph is simply a
graph with a mapping of nodes to points in $\R^D$. Based on this approximation,
we re-derive the asymmetric weights for lattices described in the
previous section as values of continuous functions defined over
$\R^D$ along the principal axes in $\R^D$. In a lattice, the neighbors
of a node lie along the principal canonical axes of $\R^D$. For an arbitrary graph, the
weights are now chosen as samples of the same functions, along
directions in which the neighbors lie. 

The method is applicable to arbitrary dimension, but
we only consider the 2-D case in this paper. Graphs with 2-D drawings
are one of the most relevant classes of graphs for sensor networks
where consensus is likely to find application.

\begin{figure}
\begin{center}
\psfrag{x1}{\scriptsize $x_1$}
\psfrag{x2}{\scriptsize $x_2$}
\psfrag{o}{\scriptsize$o$}
\psfrag{1}{\tiny $1$}  
\psfrag{L}{\scriptsize $\L$}  
\includegraphics[scale = 0.3]{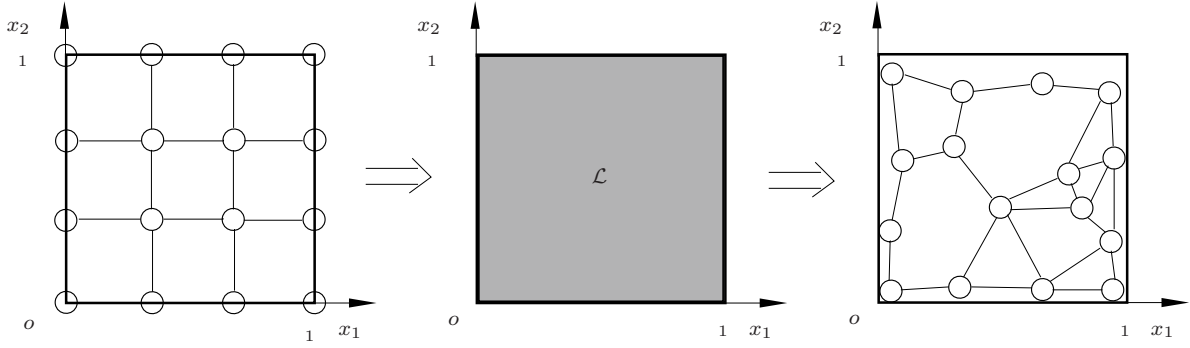} 
\end{center}
\caption{Continuum approximation of general graphs.}\label{fig:approximation}
\vspace{-5 mm}
\end{figure}

%%%%%%%%%%%%%%%%%%%%%%%%%%%%%%%%%%%%%%%%%%%%%%%%%%%%%%%%%%%%%%%%%%%%%
\subsection{Continuum approximation}
Recall that  the convergence rate is intimately connected to the Laplacian matrix. 
We will show that the Laplacian matrix associated with a large 2-D lattice  with certain weights
 can be approximated by a Sturm-Liouville operator defined on a 2-D plane. Thus
it's reasonable to suppose that the Sturm-Liouville operator
 is also a good (continuum) approximation of the Laplacian
matrix of  large graphs with 2-D drawing. We start from 2-D lattice graph and derive a Sturm-Liouville operator. We then use  this operator to approximate the graph Laplacian of more general graphs. The idea is illustrated in Figure~\ref{fig:approximation}. 

For ease of description, we first consider a 1-D lattice, with the following asymmetric weights  inspired by~\cite{HH_PB_CDC:11}, 
\begin{align}\label{eq:weight_1d}
	W_{i,i+1}=c=\frac{1+\varepsilon}{2}, \quad W_{i+1,i}=a=\frac{1-\varepsilon}{2},
\end{align}
where $i\in \{1,2,\cdots, N-1\}$ and $\varepsilon \in (0,1)$ is a constant. The graph Laplacian corresponding to the weights given in~\eqref{eq:weight_1d} is given by
\begin{align}\label{eq:laplacian_1D_v}
	L^{(1)}=\begin{bmatrix}
		\frac{1+\varepsilon}{2} & \frac{-1-\varepsilon}{2} & & \\
		\frac{-1+\varepsilon}{2} & 1 & \frac{-1-\varepsilon}{2} &\\
	        &\ddots	&\ddots	&\ddots	\\
		& & \frac{-1+\varepsilon}{2} & 1 & \frac{-1-\varepsilon}{2} \\
		& &  & \frac{-1+\varepsilon}{2} & \frac{1-\varepsilon}{2} \\
	\end{bmatrix}.
\end{align}
Recall that to find a lower bound of the convergence rate of the weight matrix $W^{(1)}$, it's sufficient to find a lower bound of the second smallest eigenvalue of the associate Laplacian matrix $L^{(1)}$.

We now use a Sturm-Liouville operator to approximate the Laplacian matrix $L^{(1)}$. We first consider the finite-dimensional eigenvalue problem $L^{(1)} \phi=\mu \phi$,
\begin{align*}
	 \begin{bmatrix}
		\frac{1+\varepsilon}{2} & \frac{-1-\varepsilon}{2} & & \\
		\frac{-1+\varepsilon}{2} & 1 & \frac{-1-\varepsilon}{2} &\\
	        &\ddots	&\ddots	&\ddots	\\
		& & \frac{-1+\varepsilon}{2} & 1 & \frac{-1-\varepsilon}{2} \\
		& &  & \frac{-1+\varepsilon}{2} & \frac{1-\varepsilon}{2} \\
	\end{bmatrix} 
	\begin{bmatrix}
	\phi_1 \\
	\phi_2\\
	\vdots	\\
	\phi_{N-1}\\
	\phi_N\\
	\end{bmatrix} =\mu 	
	\begin{bmatrix}
	\phi_1 \\
	\phi_2\\
	\vdots	\\
	\phi_{N-1}\\
	\phi_N\\
	\end{bmatrix}.
\end{align*}

Expanding the equation, we have the following coupled difference equations
\begin{align*}
	\frac{-1+\varepsilon}{2} \phi_{i-1} + \phi_i +\frac{-1-\varepsilon}{2} \phi_{i+1}= \mu \phi_i,
\end{align*}
where $i \in \{1,2,\cdots, N\}$ and $\phi_0=\phi_1$, $\phi_{N+1}=\phi_N$. The above equation can be rewritten as
\begin{align*}
	-\frac{1}{2N^2} \frac{\phi_{i-1}-2 \phi_i+\phi_{i+1}}{1/N^2} -\frac{\varepsilon}{N} \frac{\phi_{i+1}- \phi_{i-1}}{2/N}= \mu \phi_i.
\end{align*}

The starting point for the continuum approximation is to consider a function
$\phi(x): [0 , 1]  \to \R$ that satisfies:
\begin{align}\label{eq:p_approx}
  \phi_{i} = \phi(x)|_{x = i/(N+1)}, 
\end{align}
such that functions that are defined at discrete points $i$ will be approximated by functions that are defined everywhere in $[0,1]$. The original functions are thought of as samples of their continuous approximations. Under the assumption that $N$ is large, using the following finite difference approximations:
\begin{align*}
\Big [ \frac{\phi_{i-1}-2\phi_{i}+\phi_{i+1}}{1/N^2} \Big ]
&=\Big [\frac{\partial^2 \phi(x,t)}{\partial {x}^2} \Big ]_{x = i/(N+1)}, \\
\Big [ \frac{\phi_{i+1}-\phi_{i-1}}{2/N} \Big ]
&=\Big [\frac{\partial \phi(x,t)}{\partial {x}} \Big ]_{x =i/(N+1)},
\end{align*}
the finite-dimensional eigenvalue problem can be approximated by the following Sturm-Liouville eigenvalue problem
\begin{align}\label{eq:SL}
	\L^{(1)} \phi(x)= \mu \phi(x),\quad \L^{(1)}=- \frac{1}{2N^2} \frac{d^2}{d x^2} - \frac{\varepsilon}{N} \frac{d }{d x},
\end{align}
with the following Neumann boundary conditions
\begin{align}\label{eq:BC}
	\frac{d \phi(0)}{d x}=\frac{d \phi(1)}{d x}=0.
\end{align}

\begin{lemma}\label{lem:SL}
	The  eigenvalues of the Sturm-Liouville operator $\L^{(1)}$~\eqref{eq:SL} with boundary condition~\eqref{eq:BC} for $0<\varepsilon<1$ are real and the first two smallest eigenvalues satisfy
\begin{align*}
	\mu_1 (\L^{(1)})=0, \quad \mu_2 (\L^{(1)})\geq \varepsilon^2/2. \tag*{\frqed}
\end{align*}
\end{lemma}

\smallskip
We see from Lemma~\ref{lem:SL} that the second smallest eigenvalue of the Sturm-Liouville operator $\L^{(1)}$ is uniformly bounded away from zero. This result is not surprising, since it's a continuum counterpart of Lemma~\ref{lem:eig_1}, which shows that the second smallest eigenvalue corresponding to the 1-D lattice with designed asymmetric weights is uniformly bounded below. 

We now consider the distributed consensus on D-dimensional  lattices. In particular, we consider the following weights on the graph
\begin{align}\label{eq:weights_geo}
	W^{(D)}_{i,i^{d+}}=c_d=\frac{1+\varepsilon}{2D},\quad W^{(D)}_{i,i^{d-}}=a_d=\frac{1-\varepsilon}{2D}, 
\end{align}
where $\varepsilon \in (0,1)$ is a constant.

The Laplacian matrix of the D-dimensional square lattices with the weights given in~\eqref{eq:weights_geo} is given by $L^{(D)}=I-W^{(D)}$. Following the similar procedure as the $1$-dimensional lattice, the second smallest eigenvalue of the Laplacian matrix $L^{(D)}$ can be approximated by that of the following Sturm-Liouville operator
\begin{align}\label{eq:SL_D}
	\L^{(D)}=-\sum_{\ell=1}^D ( \frac{1}{2DN_d^2} \frac{d^2}{d x_d^2} + \frac{\varepsilon}{DN_d} \frac{d }{d x_d}),
\end{align}
with the following Neumann boundary conditions
\begin{align}\label{eq:BC_D}
	\frac{\partial \phi(\vec{x})}{\partial x_d}\Big |_{x_d  = 0 \text { or } 1}=0,
\end{align}
where $d=1,2, \cdots, D$ and $\vec{x}=[x_1,x_2, \cdots, x_D]^T$.

The continuum approximation has been used to study the stability margin of large vehicular platoons~\cite{HH_PB_CDC:11,HH_PB_PM_TAC:11}, in which the continuum model gives more insight on the effect of asymmetry on the stability margin of the systems. In this paper, we use the second smallest eigenvalue of the Sturm-Liouville operator $\L^{(D)}$ to approximate that of the Laplacian matrix $L^{(D)}$.

\begin{theorem}\label{thm:SL_D}
	The second smallest eigenvalues $\mu_2 (\L^{(D)})$ of the Sturm-Liouville operator $\L^{(D)}$~\eqref{eq:SL_D} with boundary condition~\eqref{eq:BC_D} for $0<\varepsilon<1$ is real and satisfies
\begin{align}
	\mu_2 (\L^{(D)})\geq \frac{\varepsilon^2}{2D},
\end{align}
which is a positive constant independent of $N$. \frqed
\end{theorem}

%%%%%%%%%%%%%%%%%%%%%%%%%%%%%%%%%%%%%%%%%%%%%%%%%%%%%%%%%%%%%%%%%%%%%
\subsection{Weight design for general graphs}
\begin{figure}
\begin{center}
\psfrag{x1}{\scriptsize $x_1$}
\psfrag{x2}{\scriptsize $x_2$}
\psfrag{o}{\scriptsize $o$}
\psfrag{q}{\scriptsize $\theta_{12}$}
\psfrag{q1}{\scriptsize $\theta_{13}$}
\psfrag{0}{\scriptsize $0$} \psfrag{1}{\scriptsize $1$}  \psfrag{2}{\scriptsize $2$}  \psfrag{3}{\scriptsize $3$}  
\psfrag{L}{\scriptsize $\L$}  
\psfrag{11}{\scriptsize $\frac{\pi}{2}$}  
\psfrag{22}{\scriptsize $\pi$}  
\psfrag{33}{\scriptsize $\frac{3\pi}{2}$}  
\psfrag{66}{\scriptsize $2\pi$}  
\psfrag{44}{\scriptsize $\frac{1-\varepsilon}{4}$} 
\psfrag{55}{\scriptsize $\frac{1+\varepsilon}{4}$}
\psfrag{x}{\scriptsize$\theta$}
\psfrag{y}{\scriptsize $g$}

\subfigure[Relative angle]{\includegraphics[scale = 0.35]{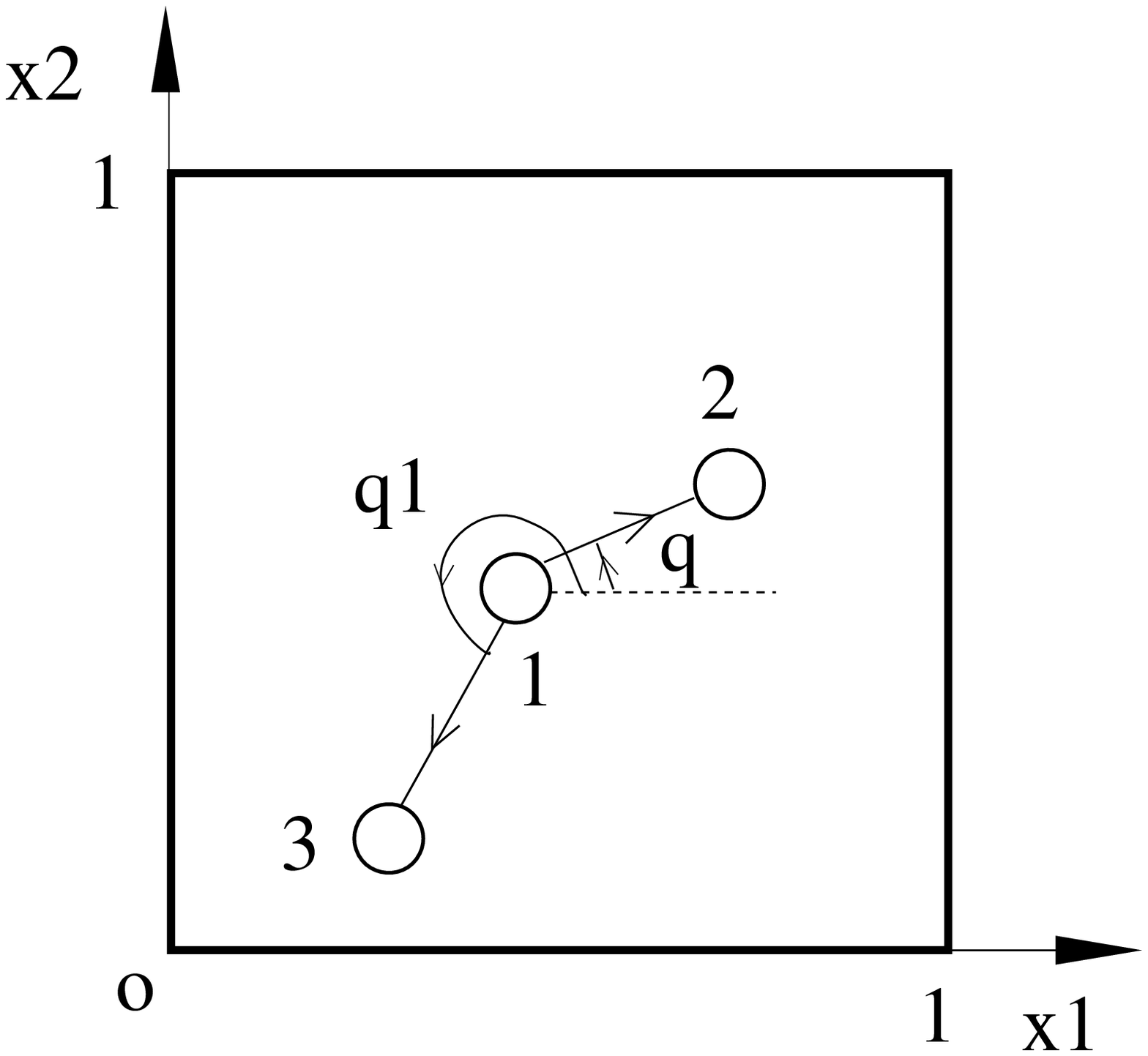}} \qquad
\subfigure[Weight function]{\includegraphics[scale = 0.35]{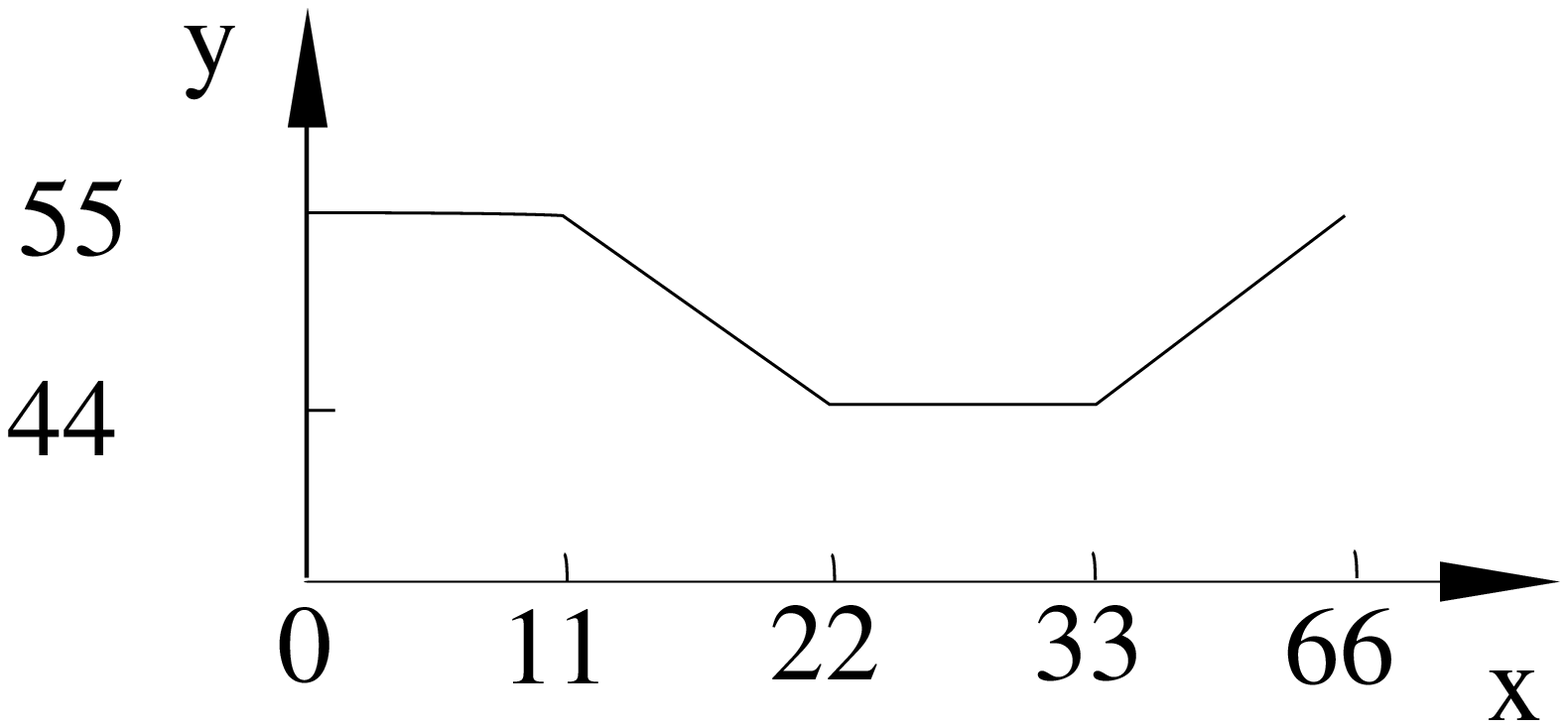}}
\end{center}
\caption{Weight design for general graphs.}\label{fig:weight}
\vspace{-5 mm}
\end{figure}

The inspiration of the proposed method comes from the design for
lattices. The $4$ weights for each node $i$ in a 2-D lattice can be
re-expressed as samples of a continuous function $g: [0, 2\pi) \to
[\frac{1-\epsilon}{4}, \frac{1+\epsilon}{4}]$:
\begin{align*}
	W_{i,i^{1+}}&= g(\theta_{i,i^{1+}}), & W_{i,i^{2+}}&=
        g(\theta_{i,i^{2+}}), \\
          W_{i,i^{1-}} &= g(\theta_{i,i^{1-}}), & W_{i,i^{2-}}&= g(\theta_{i,i^{2-}})
\end{align*}
where $\theta_{i,j}$ is the relative
angular position of $j$ with respect to $i$. Given the angular
positions of $i$'s neighbors and the values of the weights, we know
that the function $g$ must satisfy:
\begin{align}\label{eq:weight_design1}
g([0, \frac{\pi}{2}, \pi,\frac{3\pi}{2}]) &=[\frac{1+\varepsilon}{4}, \frac{1+\varepsilon}{4}, \frac{1-\varepsilon}{4}, \frac{1-\varepsilon}{4}].
\end{align}
Thus, we choose the function $g$ as shown in Figure~\ref{fig:weight} (b).  

For an arbitrary graph, we now choose the weights by sampling the
function according to the angle associated with each edge $(i,j)$:
\begin{align}\label{eq:weight_design}
	W_{i,k}=\frac{g(\theta_{i,k})}{\sum_{j \in \N_i} g(\theta_{i,j})},
\end{align}
where $g(\cdot)$ is the function described in Figure~\ref{fig:weight}
(b). The above weight function~\eqref{eq:weight_design} can be seen as
a linear interpolation of~\eqref{eq:weight_design1}. We see
from~\eqref{eq:weight_design} that the weight on each edge is
computable in a distributed manner; a node only needs to know the
angular position of its neighbors. This design method does not require
any knowledge of the network topology or centralized computation.

%%%%%%%%%%%%%%%%%%%%%%%%%%%%%%%%%%%%%%%%%%%%%%%%%%%%%%%%%%%%%%%%%%%%%
\subsection{Numerical comparison}
In this section, we present the numerical comparison of  convergence
rates among asymmetric design, symmetric optimal weights and weights
chosen by the Metropolis-Hastings method. The symmetric optimal weights
are obtained by using convex optimization
method~\cite{boyd2004fastest,xiao2004fast}. The Metropolis-Hastings
 weights are picked by the following rule: $W_{i,j}=1/|\N_i|$, where $\N_i$ denotes the number of node $i$'s neighbors. The weights generated by this method are in general asymmetric. We plot the convergence rate $R$ as a function of $N$, where $N$ is the number of agents in the network. The amount of asymmetry used is $\varepsilon=0.5$.

\begin{figure}[h]
	  \psfrag{R}{\scriptsize $R$}
	  \psfrag{N}{\scriptsize $N$}
	  \psfrag{Convex optimization}{\scriptsize Symmetric optimal}
	  \psfrag{Asymmetric design}{\scriptsize Asymmetric Design}
	  \psfrag{Metropolis-Hastings}{\scriptsize Metropolis-Hastings}
	  \psfrag{0}{ } \psfrag{0.5}{ } \psfrag{1}{ }
\centering
\subfigure[L-Z geometric]{\includegraphics[scale = 0.27]{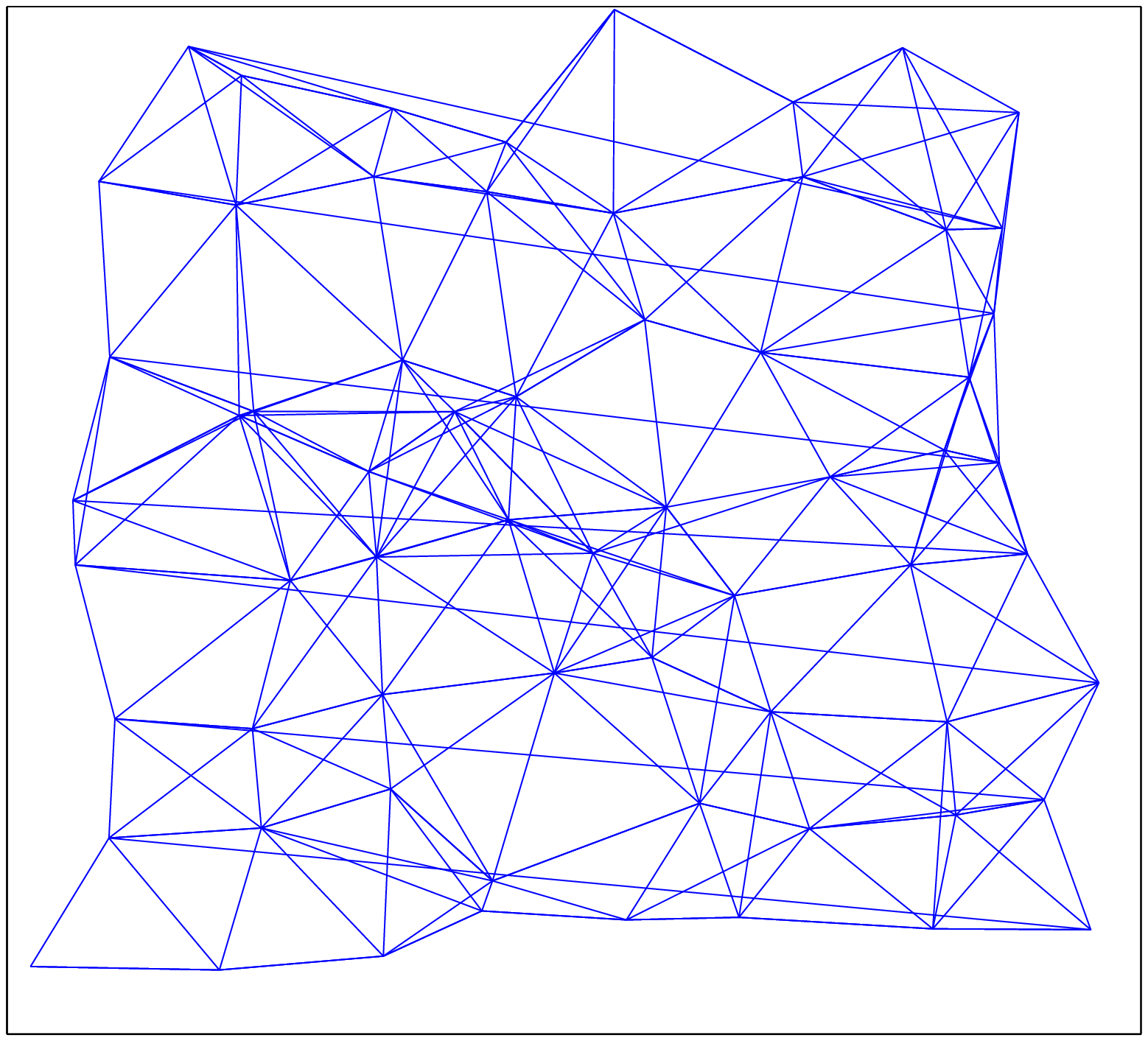}}
\subfigure[Delaunay]{\includegraphics[scale = 0.27]{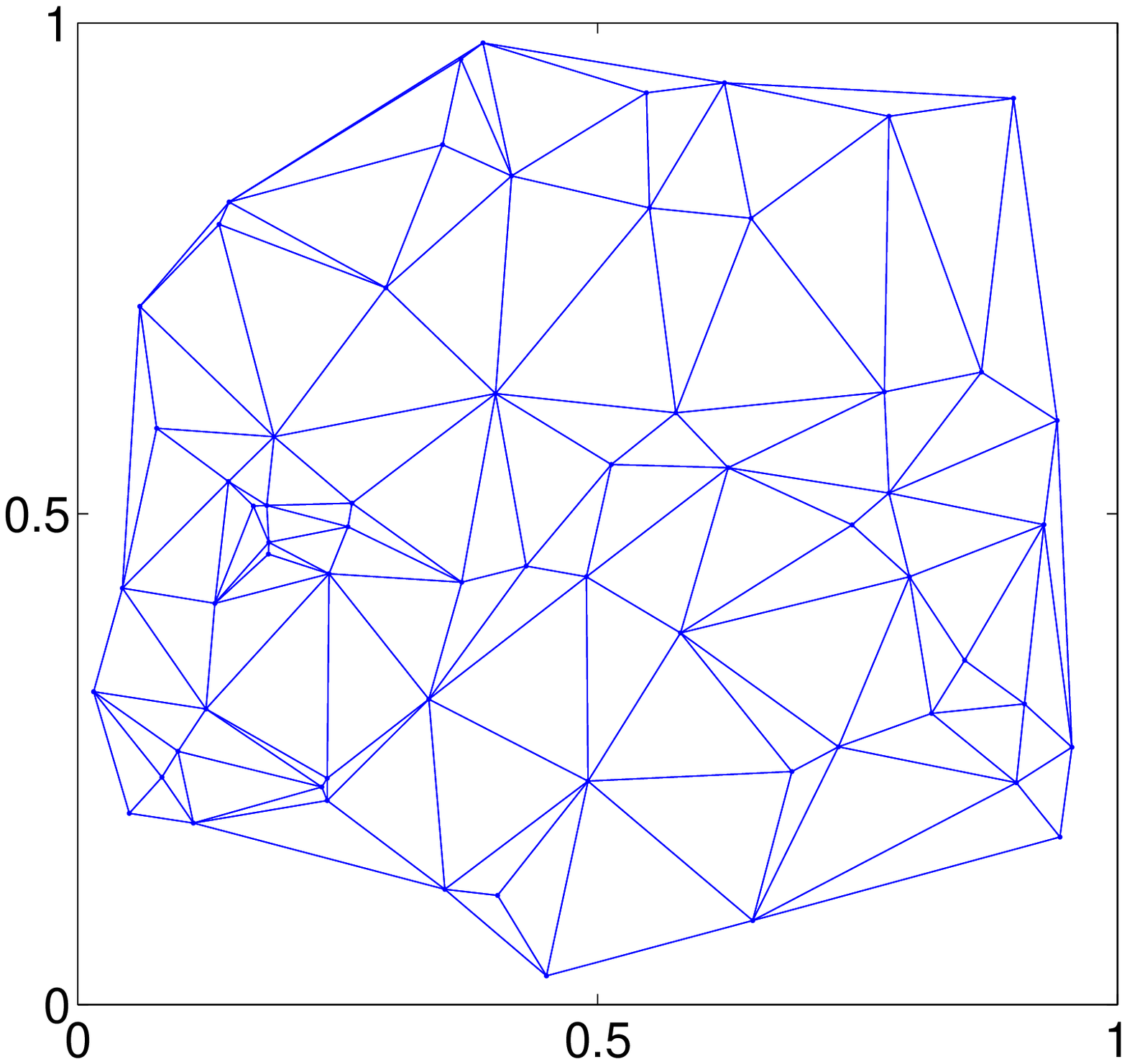}}
\subfigure[Random geometric]{\includegraphics[scale = 0.27]{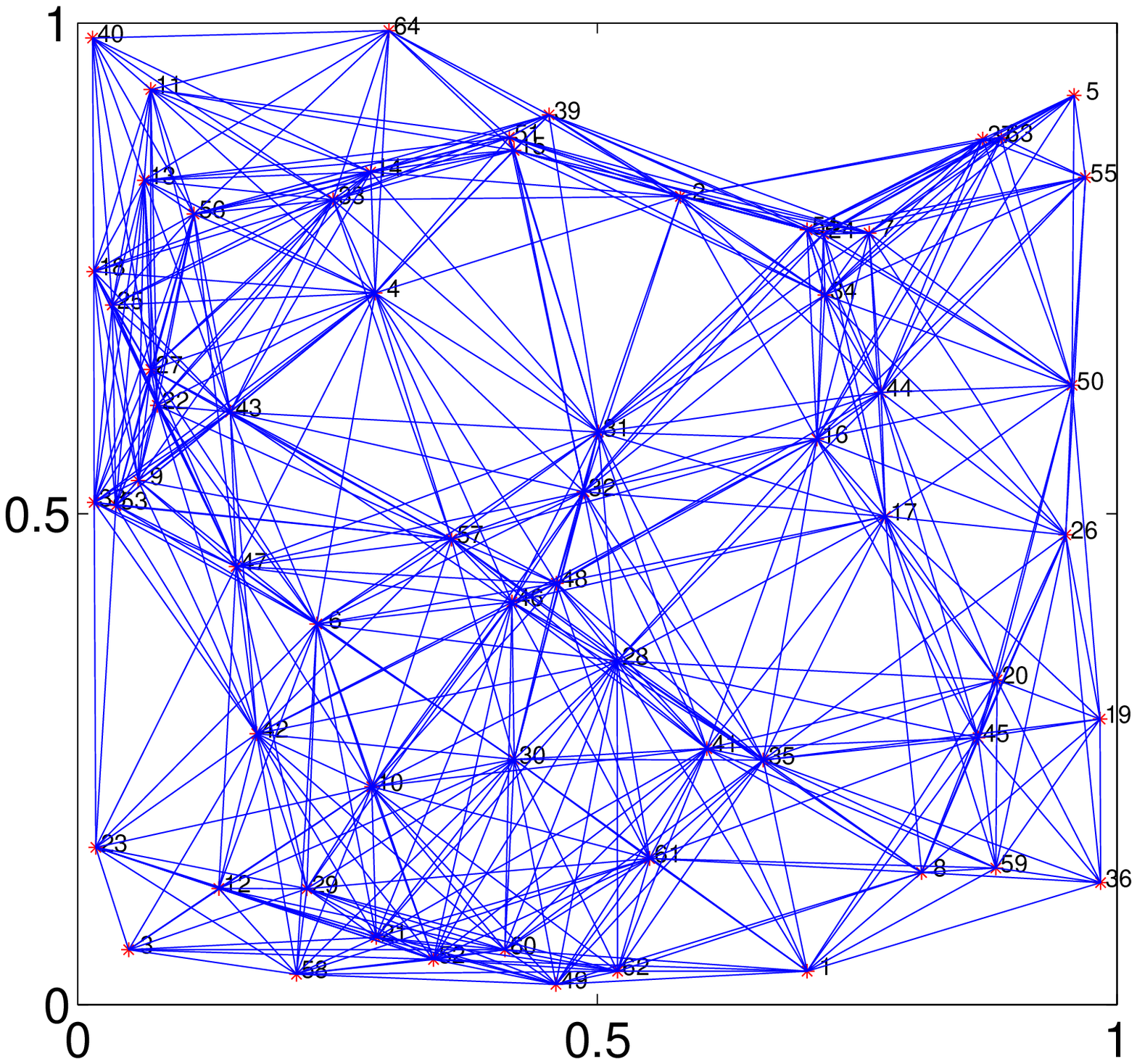}}
\caption{Examples  of $2$-D L-Z geometric, Delaunay and random geometric graphs.}
\label{fig:ex}
%\vspace{-5 mm}
\end{figure}

We first consider a L-Z geometric graph~\cite{LovisariARC2011}, which
is generated by perturbing the position of a square $2$-D lattice
($N_1=N_2=\sqrt{N}$) with Gaussian random noise (zero mean and
$1/(4\sqrt{N})$ standard deviation) and connect each nodes with the
other nodes that are within $2/\sqrt{N}$ of radius
neighborhood. Second, we consider a Delaunay
graph~\cite{martinez2007motion}, which is generated by placing $N$
nodes on a 2-D unit square uniformly at random and connecting any two
nodes if their corresponding Voronoi cells intersect, as long as their
Euclidean distance is smaller than $1/3$. Finally, we consider a random
geometric graphs~\cite{boyd2006randomized}, which is generated by
placing $N$ nodes on a 2-D unit square uniformly at random and
connecting pairs of nodes that are within distance $3/\sqrt{N}$ of each
other. Figure~\ref{fig:ex} gives examples of L-Z geometric graphs,
Delaunay graphs and random geometric graphs.

\begin{figure}
	  \psfrag{R}{\scriptsize $R$}
	  \psfrag{N}{\scriptsize $N$}
	  \psfrag{Convex optimization}{\scriptsize Symmetric optimal}
	  \psfrag{Asymmetric design}{\scriptsize Asymmetric Design}
	  \psfrag{Metropolis-Hastings}{\scriptsize Metropolis-Hastings}
\centering
\subfigure[L-Z geometric graphs]{\includegraphics[scale = 0.35]{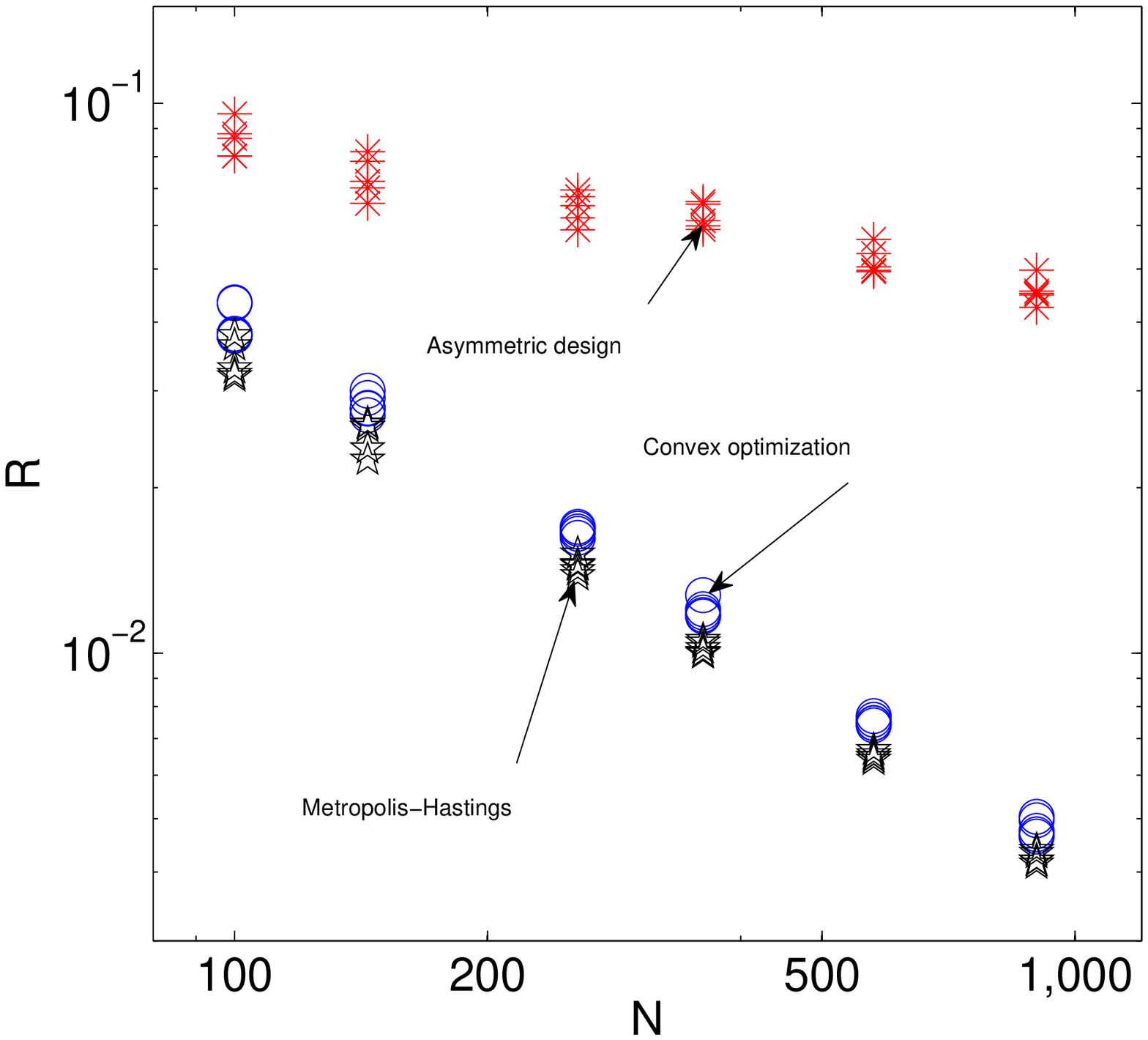}} \\
\subfigure[Delaunay graphs]{\includegraphics[scale = 0.35]{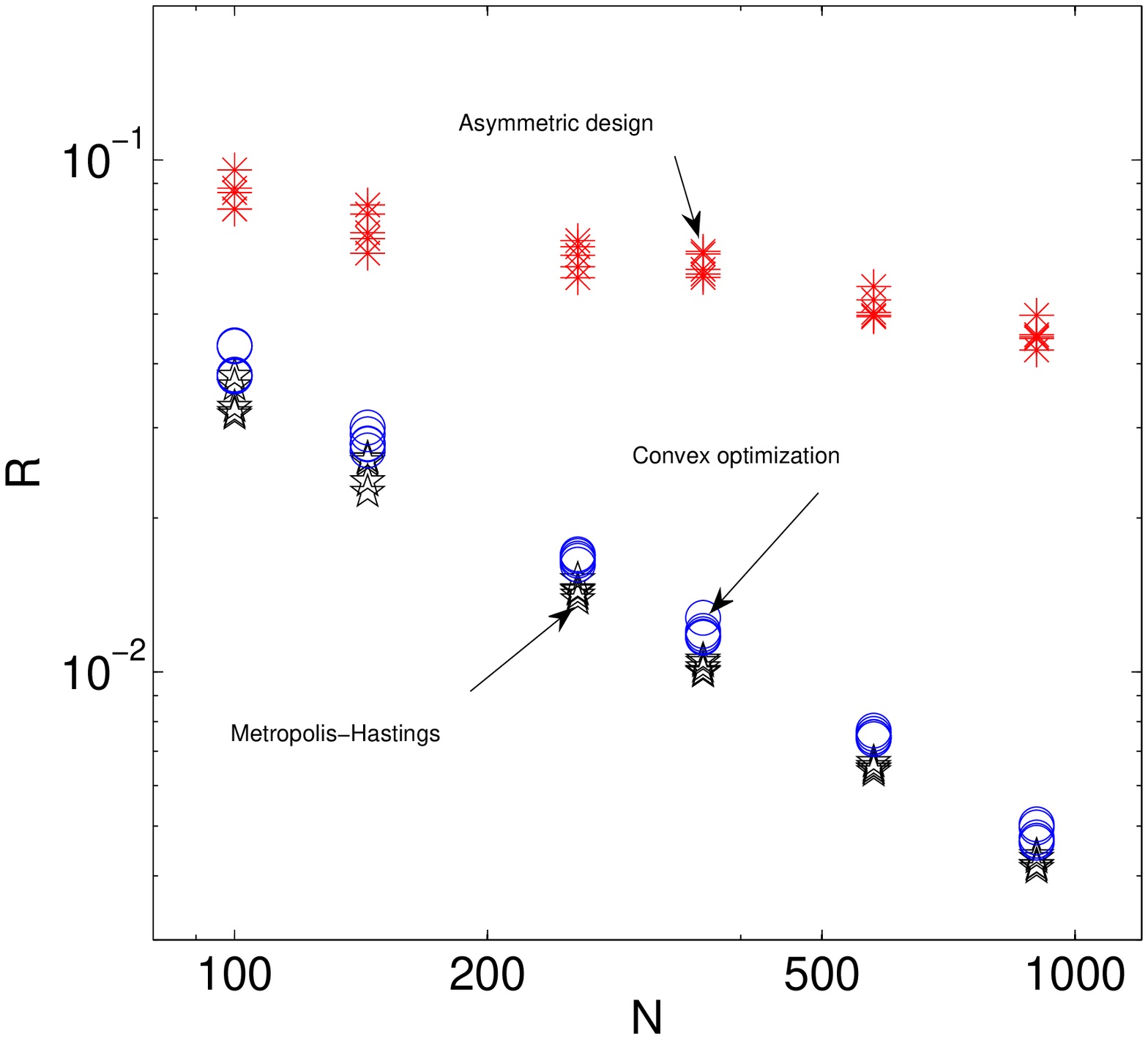}} \\
\subfigure[Random geometric graphs]{\includegraphics[scale = 0.35]{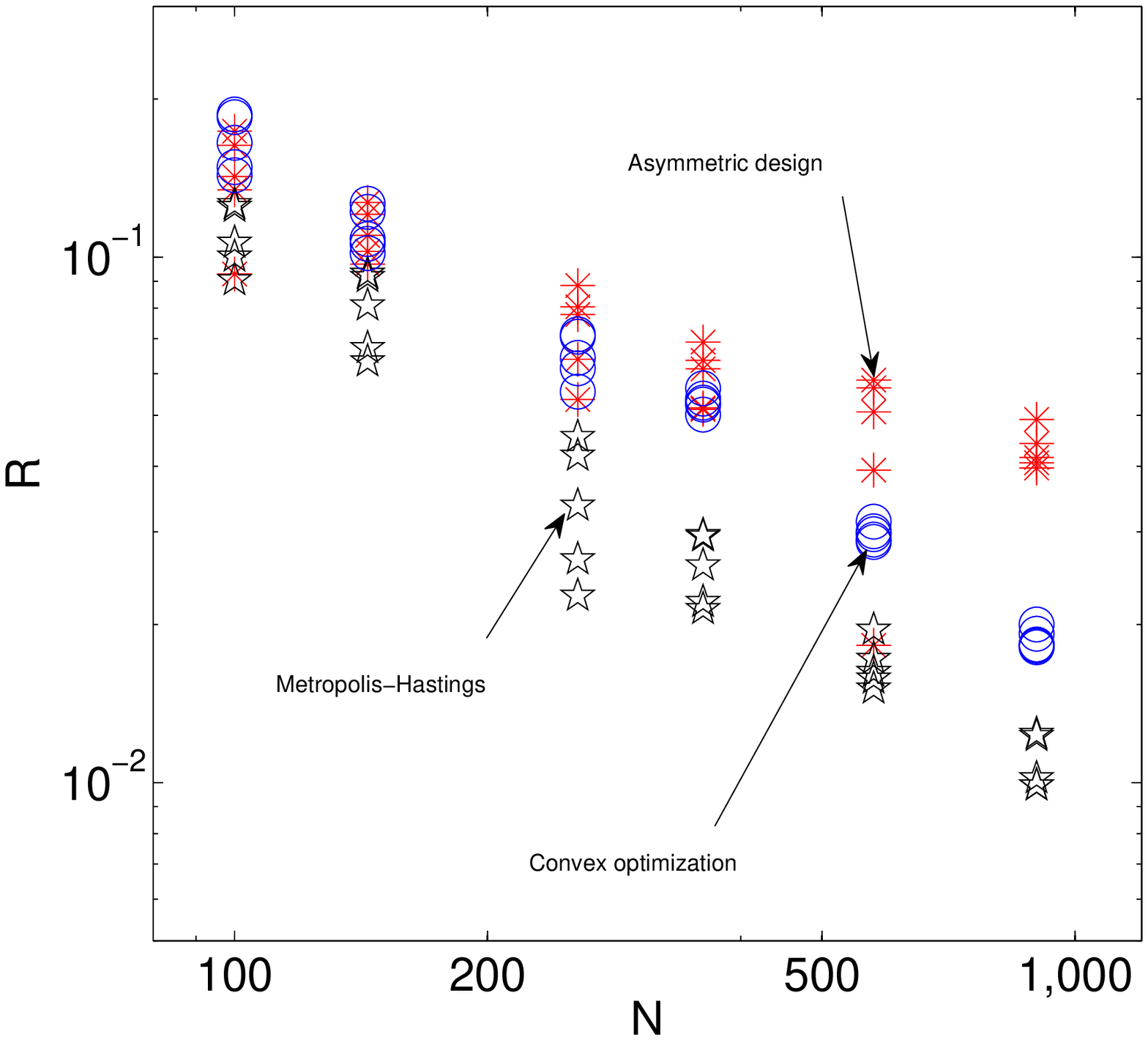}}

\caption{Comparison of convergence rates with proposed asymmetric
  weights, Metropolis-Hastings weights, and symmetric optimal. For
  each $N$, results from $5$ sample graphs are plotted.}\label{fig:comparison}
\end{figure}

\begin{figure}
 	  \psfrag{R}{\scriptsize $E[R]$}
 	  \psfrag{N}{\scriptsize $N$}
 	  \psfrag{AS}{\scriptsize Asymmetric design}
 	  \psfrag{MH}{\scriptsize Metropolis-Hastings}
 \centering
 \subfigure[L-Z geometric]{\includegraphics[scale = 0.265]{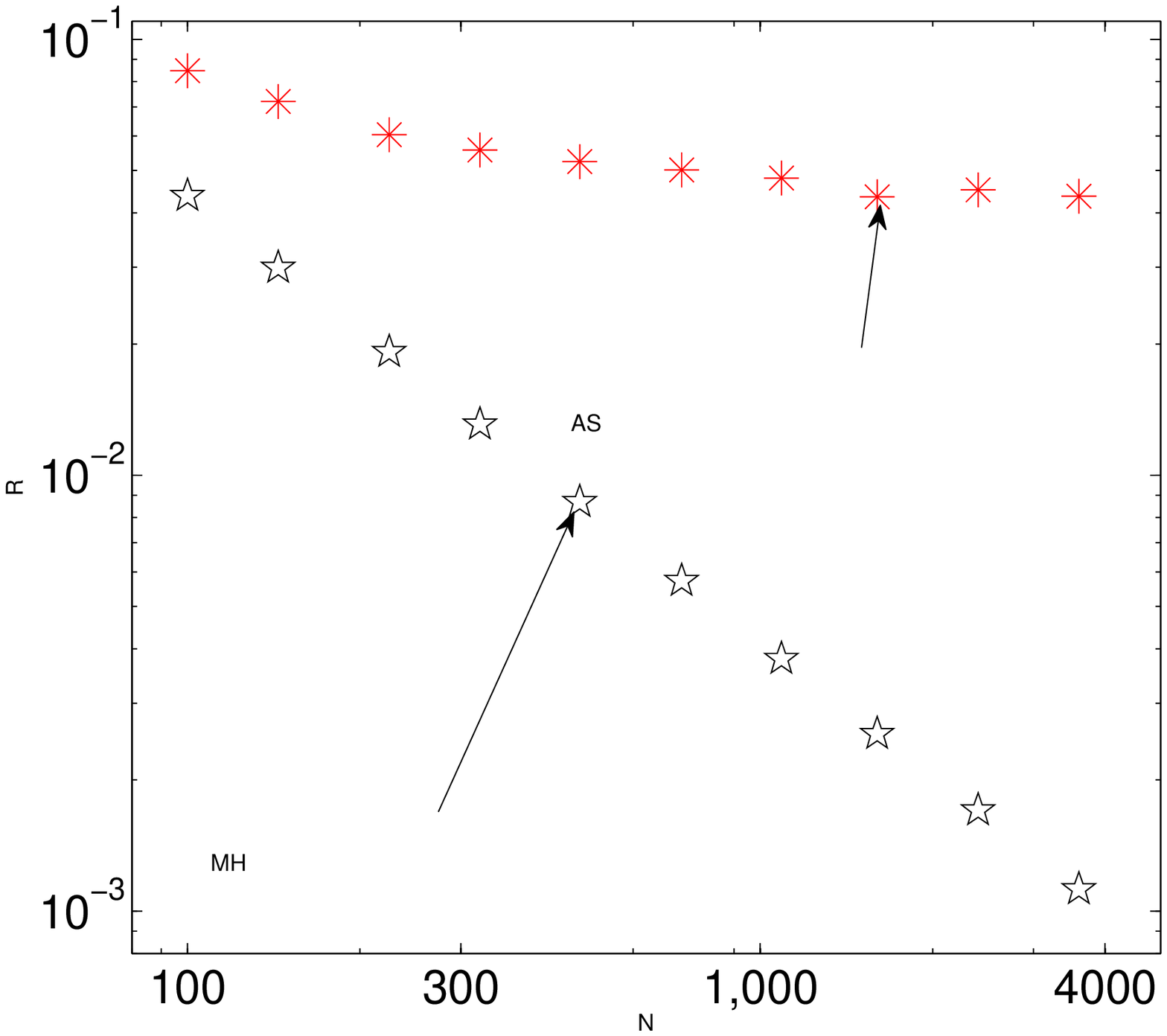}}
 \subfigure[Delaunay]{\includegraphics[scale = 0.265]{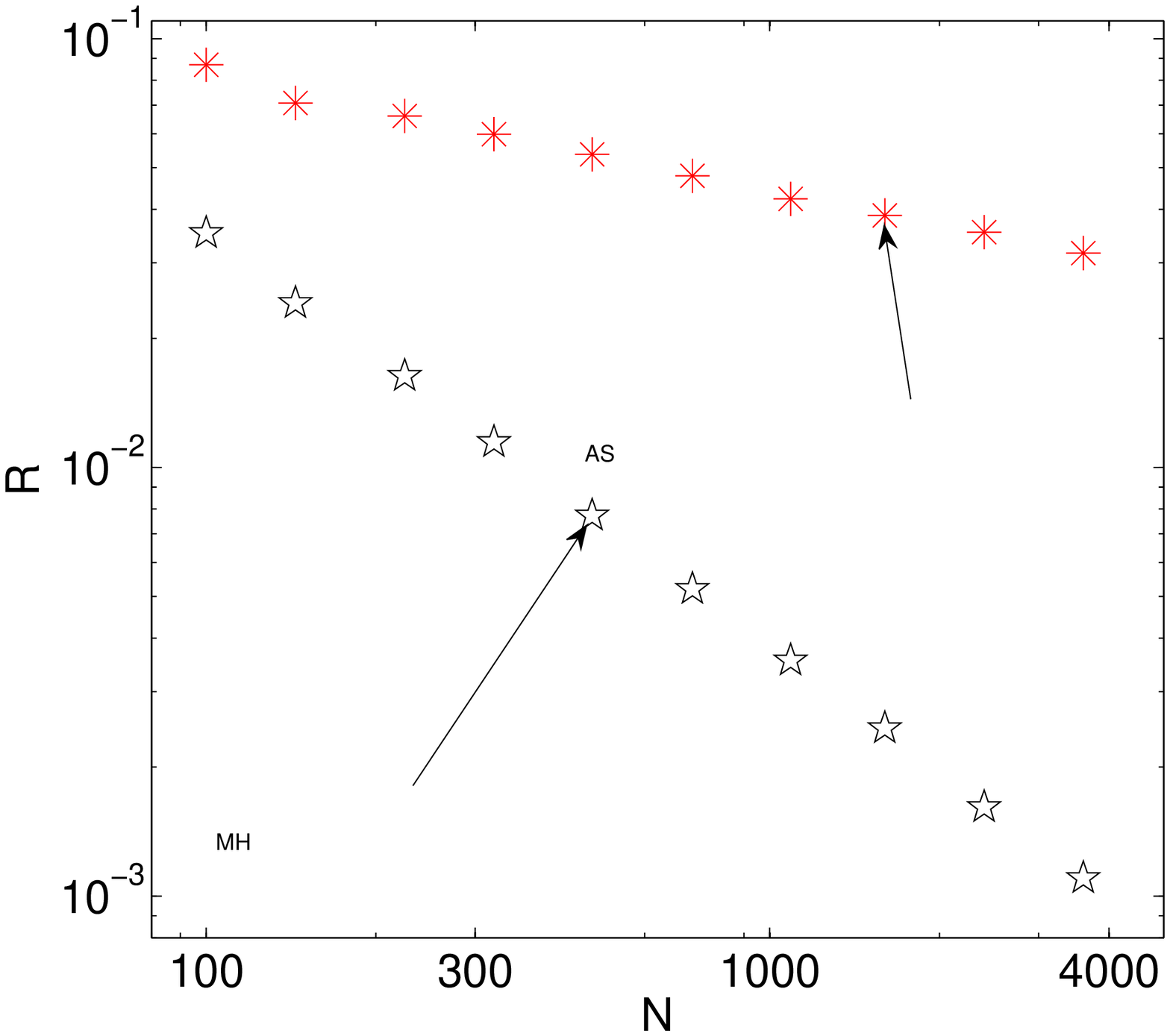}}
 \subfigure[Random geometric]{\includegraphics[scale = 0.265]{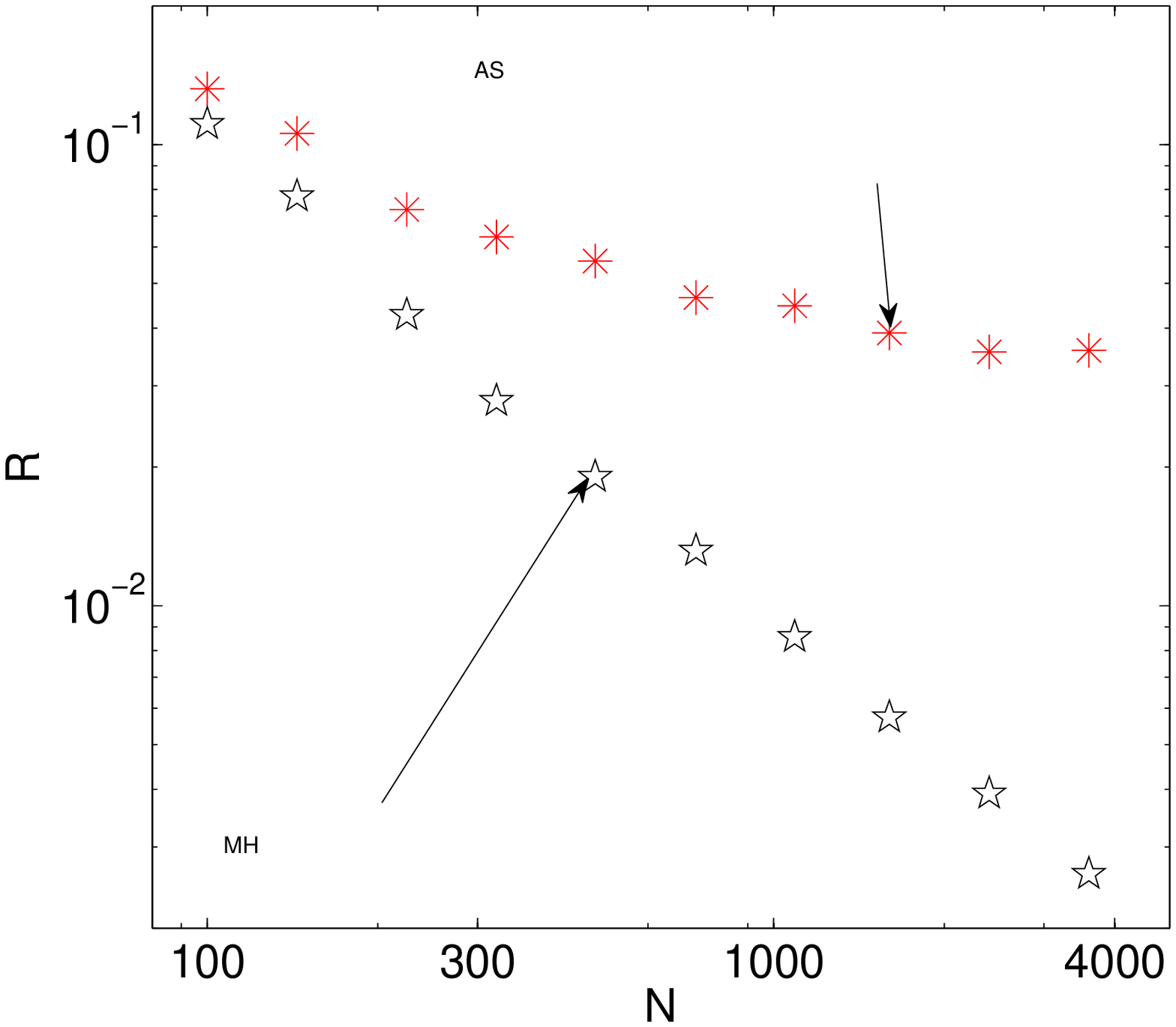}}

 \caption{Mean of convergence rates of L-Z geometric, Delaunay and random geometric graphs with asymmetric design (AS) and Metropolis-Hastings (MH) for large $N$ computed from $10$ samples.}

 \label{fig:trend}
 \end{figure}

Figure~\ref{fig:comparison} shows the comparison of convergence rates
among asymmetric design, symmetric optimal and Metropolis-Hastings
weights. For each $N$, the convergence rate of  $10$ samples of the
graphs are plotted. We see from Figure~\ref{fig:comparison}
that for almost every sample in each of the three classes, the
convergence rate with the asymmetric design is an order of magnitude
larger than the others, especially when $N$ is large.  In addition, the
 convergence rates with symmetric optimal and Metropolis-Hastings
 methods are similar.  Moreover, we
 observe from Figure~\ref{fig:trend} (a) and (c) that the slopes of the
 convergence rates with asymmetric design for L-Z geometric graphs and
 random geometric graphs are becoming progressively smaller with
 increasing $N$, which indicates that the convergence rate has a potential to be  
 uniformly bounded below when $N$ becomes arbitrarily  large.
 The convergence rates with  symmetric optimal weights are not included, since the numerical computations for large $N$ are extremely expensive.

\section{Conclusions and Future work}\label{sec:conc}
We studied the problem of how to design weights to increase the
convergence rate of distributed consensus in networks with static
topology. We proved that on lattice graphs, with proper choice of
asymmetric weights, the convergence rate can be uniformly bounded away
from zero. In addition, we propose a distributed weight design
algorithm for 2-dimensional geometric graphs  to
improve the convergence rate, by using a continuum
approximation. Numerical calculations show that the resulting
convergence rate is substantially larger than that optimal symmetric
weights and Metropolis Hastings weights.

An important open question is a precise
characterization of graphs for which theoretical guarantees on
size-independent convergence rate can be provided with the proposed
design. In addition, characterizing the asymptotic steady state value
for more general graphs than lattices is also on-going work.

\bibliographystyle{IEEEtran}        % Include this if you use bibtex 

%HH
\bibliography{../../PBbib/vehicular_platoon,../../PBbib/sensnet_bib_dbase,../../PBbib/Barooah,../../PBbib/HH,../../PBbib/distributed_control}

\appendix
\begin{proof-lemma}{\ref{lem:eig_1}}
	The stochastic matrix $W^{(1)}$ has a simple eigenvalue  $\lambda_1=1$. Following Theorem 3.1 of~\cite{yueh2008explicit}, the other eigenvalues of $W^{(1)}$ are given by
\begin{align*}
	 \lambda_\ell&=1-a-c+2 \sqrt{ac} \cos \theta_\ell, \quad \ell \in\{2,\cdots, N\},
\end{align*}
where $\theta_\ell$ ($\theta \neq m\pi, m \in \mathbb{Z}$, $\mathbb{Z}$ being
the set of integers) is the  root of the following equation
\begin{align*}
	 2 \sin(N\theta) cos(\theta) =(a+c) \sqrt{\frac{1}{ac}} \sin N\theta,
\end{align*}
which implies
\begin{align*}
	\sin(N\theta)=0, \quad \text{or} \quad \cos \theta = \frac{(a+c)}{2} \sqrt{\frac{1}{ac}}.
\end{align*}
Since $a>0, c>0$ and $a \neq c$, we have $\frac{(a+c)}{2} \sqrt{\frac{1}{ac}} >1$, thus $\cos \theta \neq \frac{(a+c)}{2} \sqrt{\frac{1}{ac}}$. In addition, we have that $\theta \neq m\pi$, which yields
\begin{align}
\theta_\ell =\frac{(\ell-1) \pi}{N}, \ \ell=\{2, \cdots, N\}.
\end{align}
We now obtain the eigenvalues of $W^{(1)}$, which is given by
\begin{align*}
	\lambda_\ell =1-a-c+2\sqrt{ac} \cos \frac{(\ell-1) \pi}{N}, \ \ell=\{2, \cdots, N\}. 
\end{align*}

 Let $\pi=[\pi_1, \pi_2, \cdots, \pi_N]$ be the left Perron vector of $W^{(1)}$. From the definition of Perron vector, we have $\pi W^{(1)}= \pi$.  Thanks to the special structure of the tridiagonal form of $W^{(1)}$, we can solve for $\pi$ explicitly, which yields
\begin{align}\label{eq:perron_element}
	\pi_i=(c/a)^{i-1} \pi_1,
\end{align}
where $i \in \{2, 3, \cdots, N\}$.
In addition, we have $\pi_i>0$ and $\sum_{i=1}^N \pi_i=1$. Therefore,
\begin{align*}
	1= \sum_{i=1}^N \pi_i =  \sum_{i=1}^N (c/a)^{i-1} \pi_1 \quad \Rightarrow \quad \pi_1 = 	\frac{1-c/a}{1-(c/a)^N}.
\end{align*}
Substituting the above equation into~\eqref{eq:perron_element}, we complete the proof. 	\frQED
\end{proof-lemma}

\begin{proof-lemma}{\ref{lem:eig_D}}
With the weights given in~\eqref{eq:weights}, it is straightforward - through a bit tedious - to show that the graph Laplacian $L^{(D)}$ associated with the $D$-dimensional lattice  has the following form:
\begin{align*}
	L^{(d)} = I_{N_d} \otimes L^{(d-1)}+ L^{(1)}_d \otimes I_{N_1 N_2\cdots N_{d-1}}, \quad 2 \leq d \leq D,
\end{align*}
where $L^{(1)}=L^{(1)}_1$ and $L^{(1)}_d=1-W^{(1)}_d$ is the Laplacian matrix of dimension $N_d \times N_d$, which is given by
\begin{align}\label{eq:triangular3}
	L^{(1)}_d=\begin{bmatrix}
		c_d & -c_d & & \\
		-a_d & a_d+c_d & -c_d &\\
	        &\ddots	&\ddots	&\ddots	\\
		& & -a_d & a_d+c_d & -c_d \\
		& &  & -a_d & a_d \\
	\end{bmatrix}.
\end{align}

Since a $D$-dimensional lattice is the Cartesian product graph of $D$ $1$-dimensional lattices, the eigenvalues of the graph Laplacian matrix $L^{(D)}$ are sum of the eigenvalues of the $D$ $1$-dimensional Laplacian matrix $L^{(1)}_d$. Thus, we have
\begin{align*}
	\mu_{\ell_1, \dots, \ell_{D}} (L^{(D)}) =\sum_{d=1}^{D} \mu_{\ell_d} (L^{(1)}_d).
\end{align*}
In addition, we have that $W^{(D)}=I_N- L^{(D)}$ and $W^{(1)}_d=I_{N_d}- L^{(1)}_d$, thus the eigenvalues $\lambda_{\vec{\ell}}$ of $W^{(D)}$ are given by
\begin{align*}
	\lambda_{\vec{\ell}}\ (W^{(D)})&=1-\mu_{\vec{\ell}} \ (L^{(D)}) =1- \sum_{d=1}^{D} \mu_{\ell_d} (L^{(1)}_d)\\&=1-\sum_{d=1}^{D} (1-\lambda_{\ell_d} (W^{(1)}_d)).
\end{align*}

To see $\pi= \pi^{(1)}_D\otimes \pi^{(1)}_{D-1} \otimes \cdots \otimes \pi^{(1)}_1$ is the left Perron vector of $W^{(D)}$, we first notice that
\begin{align*}
	\pi^{(1)}_d W^{(1)}_d=	\pi^{(1)}_d, \quad \pi^{(1)}_d L^{(1)}_d= 0, 
\end{align*}
where $d \in \{1,\cdots, D\}$. The rest of the proof follows by straightforward induction method, we omit the proof due to space limit. \frQED
\end{proof-lemma}

\begin{proof-lemma}{~\ref{lem:SL}}
 Multiply both sides of~\eqref{eq:SL} by $2N^2e^{2\varepsilon N x}$, we obtain the  standard Sturm-Liouville eigenvalue problem
\begin{align}\label{eq:SLeigen}
	\frac{d }{d x}\Big(e^{2\varepsilon N x}\frac{d \phi(x)}{d x}\Big)+ 2N^2 \mu e^{2\varepsilon N x}\phi(x)=0.
\end{align}
According to Sturm-Liouville Theory, all the eigenvalues are real, see~\cite{haberman, evans}. To solve the Sturm-Liouville eigenvalue problem~\eqref{eq:SL}-\eqref{eq:BC}, we assume solution of the form, $\phi(x)=e^{rx}$, then we obtain the following equation
\begin{align}\label{eq:chara_1}
	r^2+2\varepsilon N r+ 2\mu N^2=0, \notag \\
	\Rightarrow    r=N(-\varepsilon \pm \sqrt{\varepsilon^2-2\mu}).
\end{align}
Depending on the discriminant in the above equation, there are three cases to analyze:
\begin{enumerate}
\item $\mu < \varepsilon^2/2$, then the eigenfunction $\phi(x)$ has the following form
		$\phi(x)=c_1 e^{N(-\varepsilon + \sqrt{\varepsilon^2-2\mu})x}+c_2 e^{N(-\varepsilon - \sqrt{\varepsilon^2-2\mu})x}$, where $c_1, c_2$ are some constants. Applying the boundary condition~\eqref{eq:BC}, it's straightforward to see that, for non-trivial eigenfunctions $\phi(x)$ to exit,  the following equation must be satisfied $$\frac{-\varepsilon  +  \sqrt{\varepsilon^2-2\mu}}{\varepsilon  +  \sqrt{\varepsilon^2-2\mu}}=e^{2N \sqrt{\varepsilon^2-2\mu}} \frac{-\varepsilon  +  \sqrt{\varepsilon^2-2\mu}}{\varepsilon  +  \sqrt{\varepsilon^2-2\mu}}.$$ Thus, we have $\mu=0$.
\item $\mu =\varepsilon^2/2$, then the eigenfunction $\phi(x)$ has the following form
	\begin{align*}
		\phi(x)=c_1 e^{-\varepsilon N x}+c_2 x e^{-\varepsilon N x}.
	\end{align*}
Applying the boundary condition~\eqref{eq:BC} again, it's straightforward to see that there is no eigenvalue for this case.
\item $\mu >\varepsilon^2/2$, then the eigenfunction has the following form
$\phi(x)= e^{-\varepsilon N x} (c_1 \cos (N\sqrt{2\mu-\varepsilon^2} x)+c_2 \sin(N\sqrt{2\mu-\varepsilon^2} x)$. 
Applying the boundary condition~\eqref{eq:BC}, for non-trivial eigenfunctions  to exit, the eigenvalues $\mu$ must
satisfy $\mu=\frac{\varepsilon^2}{2} + \frac{\ell^2 \pi^2}{2N^2}$, where $\ell = 1,2, \cdots$. 
\end{enumerate} 
Combining the above three cases, the eigenvalues of the Sturm-Liouville operator are $\mu \in \{0, \frac{\varepsilon^2}{2} + \frac{\ell^2 \pi^2}{2N^2}\}$, where $\ell \in \{1,2, \cdots\}$. The second smallest eigenvalue $\mu_2(\L)$ of the Strum-Liouville operator $\L$ is then given by
\begin{align*}
	\mu_2(\L)=\frac{\varepsilon^2}{2}+\frac{\pi^2}{2N^2}\geq \frac{\varepsilon^2}{2},
\end{align*}
which is a constant that is bounded away from $0$. \frQED
\end{proof-lemma}

\begin{proof-theorem}{\ref{thm:bound}}
According to Lemma~\ref{lem:eig_1}, the eigenvalues of $W^{(1)}_d$ are given by:
\begin{align*}
\lambda_1(W^{(1)}_d)&=1, \\ \lambda_{\ell}(W^{(1)}_d)&=1-a_d-c_d+2\sqrt{a_dc_d} \cos \frac{(\ell_d-1) \pi}{N_d}.
\end{align*}
From Lemma~\ref{lem:eig_D}, the second largest eigenvalue $\lambda_2 (W^{(D)})$ and the smallest eigenvalue $\lambda_N (W^{(D)})$ of $W^{(D)}$ are given by  
\begin{align}
	\lambda_2 (W^{(D)}) &= 1- \max_{d \in \{1,\cdots, D\}} (1-\lambda_2 (W^{(1)}_d)) \notag\\ &\leq 1- \max_{d \in \{1,\cdots, D\}} (a_d+c_d-2\sqrt{a_d c_d}), \label{eq:bound_mu2}\\ 
	\lambda_N (W^{(D)}) &= 1-\sum_{d=1}^{D} (1-\lambda_{N_d} (W^{(1)}_d)) \notag \\&=1-\sum_{d=1}^{D} (a_d+c_d-2\sqrt{a_d c_d}\cos \frac{(N_d-1) \pi}{N_d}) \notag\\ &\geq 1-\sum_{d=1}^{D} (a_d+c_d-2\sqrt{a_d c_d}).\label{eq:bound_muN}
\end{align}
Recall that $R=\min\{1-\lambda_2, 1+\lambda_N \}$. In addition, $a_d, c_d$ are fixed constants and satisfy $a_d \neq c_d$, $\sum_{d=1}^D a_d+c_d \leq 1$, therefore the lower bounds of $1-\lambda_2 (W^{(D)})$ and $1+\lambda_N (W^{(D)})$ are fixed positive constants.  We then have that the convergence rate of $W^{(D)}$ satisfy $R=1-\rho(W^{(D)}) \geq  c_0$, where $c_0$ is a constant independent of $N$.\frQED
\end{proof-theorem}

\begin{proof-theorem}{\ref{thm:SL_D}}
	By the method of separation of variables~\cite{haberman, evans}, the eigenvalues of the Sturm-Liouville operator $\L^{(D)}$ is given by
	\begin{align}\label{eq:eig_D}
	\mu (\L^{(D)})=\sum_{d=1}^D \mu (\L^{(1)}_d), 
\end{align}
where $\L^{(1)}_d$ is the $1$-dimensional Sturm-Liouville operator given by
\begin{align*}
	\L^{(1)}_d=-  \frac{1}{2DN_d^2} \frac{d^2}{d x_d^2} - \frac{\varepsilon}{DN_d} \frac{d }{d x_d},
\end{align*}
with  Neumann boundary conditions. Following Lemma~\ref{lem:SL}, we have that the smallest eigenvalue of $\L^{(1)}_d$ is $0$ and the second smallest eigenvalue of $\L^{(1)}_d$ is bounded below by $\L^{(1)}_d \geq \varepsilon^2/2D$. Therefore, we have from~\eqref{eq:eig_D} that the second smallest eigenvalue is 
\begin{align*}
	\mu_2 (\L^{(D)}) = \min_d \{\mu_2 (\L^{(d)})\} \geq \frac{\varepsilon^2}{2D}. \tag*{\frQED}
\end{align*}
\end{proof-theorem}

\end{document}